\documentclass{article}
\usepackage{amsfonts,amsmath,amsthm}
\usepackage[latin1]{inputenc}
\newtheorem{proposition}{Proposition}

\newtheorem{cor}{Corollary}
 
\newtheorem{theorem}{Theorem}

\theoremstyle{remark}
\newtheorem{remark}{Remark}

\let\epsilon=\varepsilon

\let\phi=\varphi
\let\al=\alpha

\let\si=\sigma  \let\ze=\zeta
  
\let\tilde=\widetilde

\newcommand{\bigpar}[1]{\bigl(#1\bigr)}

\newcommand{\field}[1]{\mathbb{#1}}
\newcommand{\R}{\field{R}}
\newcommand{\N}{\field{N}}
\newcommand{\Z}{\field{Z}}
\newcommand{\Q}{\field{Q}}

\newcommand{\T}{\field{T}}

\newcommand{\F}{{\mathcal{F}}}

\def\der^#1_#2{\frac{\partial^{#1}}{\partial {#2}^{#1}}}
\newcommand{\ind}[1]{\mathbf 1_{#1}}
\newcommand{\set}[1]{\left\{#1\right\}}

\title{Testing the finiteness of the support of a distribution: a statistical look at Tsirelson's equation}
\author{Sylvain Delattre\\ LPMA, Universit\'e Paris Diderot \\ sylvain.delattre@univ-paris-diderot.fr\and
Mathieu Rosenbaum\\ LPMA, Universit\'e Pierre et Marie Curie \\
  mathieu.rosenbaum@upmc.fr}

\date{25 February 2012}
\begin{document}
\maketitle
\begin{abstract}

\noindent We consider the following statistical problem: based on an i.i.d.\ sample of size $n$ 
of integer valued random variables with common law $\mu$, is it possible to test whether or not the support of $\mu$ is finite as $n$ goes to infinity?
This question is in particular connected to a simple case of Tsirelson's equation, for which it is natural to distinguish between two main configurations, the first one leading only to laws with finite support, and the second one including laws with infinite support. We show that it is in fact not possible to discriminate between the two situations, even using a very weak notion of statistical test. 

\end{abstract}

\section{Introduction} \label{intro}

\noindent In this paper, we consider an i.i.d.\ sample 
of integer valued random variables with common law $\mu$ and address the following problem: is it possible to build an asymptotic test for the finiteness (or for the non finiteness) of the support of 
$\mu$ as the sample size goes to infinity? This question is in fact motivated by Tsirelson's equation in discrete time, see \cite{ts1,ts2}. Let $(\nu_k)_{k\in -\N}$ be a sequence of probability laws  on the torus $\T=\R/\Z$. A process
$(\eta_k)_{k\in -\N}$ taking value in $\T$ is a weak solution of Tsirelson's equation associated to $(\nu_k)_{k\in -\N}$ if for all $k\in -\N$, the random variable $\xi_k=\eta_k-\eta_{k-1}$ has law $\nu_k$ and is independent of $\F^{\eta}_{k-1}=\si(\eta_n, n\le k-1)$.
If moreover for all $k\in -\N$, $\F^\eta_k=\F^\xi_k$, then the process $(\eta_k)_{k\in -\N}$ is said to be
a strong solution. It is proved by Yor in \cite{y} that there is always existence of a weak solution to the equation: if $U$ and $\ze_k$, $k\in -\N$, are independent random variables taking values in $\T$ such that $U$ is uniform and, for all $k$, $\ze_k$ has distribution $\nu_k$, then the process $(\eta^*_k)_{k\in -\N}$ defined by
\begin{equation} \label{unif}
\eta^*_0=U, \quad
\eta^*_{-n}=U-\zeta_0-\dots-\zeta_{-n+1} \quad \text{for all $n\ge 1$,}
\end{equation}
is a weak solution of Tsirelson's equation associated to $(\nu_k)_{k\in -\N}.$\\

\noindent In this work we consider only the simple case where $\nu_k=\nu$ for all $k$.
It can be shown that in terms of existence and unicity of strong or weak solutions to the equation, depending on the law $\nu$, three very different situations have to be considered, see \cite{y} for details:

\begin{itemize}
\item[-] \emph{Case 1}: {If $\nu$ is a Dirac measure:}\ \emph{existence of a strong solution, no uniqueness in law.}
\item[-] \emph{Case 2}: {If the support of $\nu$ contains at least two points and there exists an integer $p\ge 2$ and $x\in\T$ such that it is included in $x+\{k/p,~k=0,\ldots,p-1\}$:}\ \emph{no strong solution, no uniqueness in law.}
\item[-] \emph{Case 3}: If $\nu$ is a law which does not satisfy the conditions of Case 1 or Case 2: \emph{no strong solution, uniqueness in law.}
\end{itemize}

\noindent In \cite{yy,yy2}, Tsirelson's equation is seen as a cosmological caricature: $\eta_{k-1}$ represents the state of the universe at time $k-1$ and $\xi_k$ the action of the evolution process at time $k$. In particular, the authors raise
the issue of the ``initial magma" that is the description of the $\si$-algebra
$$ \F^\eta_{-\infty}=\bigcap_{n\ge 1} \si\bigpar{\eta_k, k\le -n}.$$
It turns out that it is necessarily a trivial $\si$-algebra in Case 3.\\

\noindent Our goal here is to give a statistical look at this cosmological caricature. More precisely, based on the historical observations of the state of the universe\\ $\eta_0,\eta_{-1},\dots,\eta_{-n}$, can we (at least partially) decide whether we are in Case 1, Case 2 or Case 3~?
Of course we can already remark that under Case 1, the $\xi_k$ are constant and under Case 2 the $\rho_k=\xi_{k}-\xi_{k-1}$ are all rational. Therefore, if we observe two 
different values for the $\xi_k$ in the sample, we can conclude that we are not in Case 1 and if we observe one irrational $\rho_k$,
we can discard Case 2. Consequently, if there are at least two different $\xi_k$ and
one irrational $\rho_k$, we can conclude that we are in Case 3. However, in practice, we have only access to rounded values so it does not really
make sense to build a procedure based on the fact that some quantity is irrational or not.\\ 

\noindent In order to rigorously answer the question, we use the notion of statistical test. We first consider an auxiliary problem, which is interesting on its own,
namely the possibility of building a consistent test for the finiteness (or for the non finiteness) of the support of an integer valued distribution $\mu$.
based on an i.i.d.\ sample $(X_1,\dots,X_n)$ with law $\mu$. We prove that designing such a test procedure is in fact impossible, even using a very weak notion of consistent test. Then we show that a corollary of this result is the impossibility to build a consistent test in order to separate on the one hand the union of Case 1 and Case 2 and on the other hand Case 3.\\

\noindent The paper is organized as follows. The framework we use in order to test for the finiteness (or for the non finiteness) of the support of a distribution based on i.i.d.\ data is detailed in Section \ref{Setting}. In Section \ref{Candidate}, we suggest an a priori natural functional in order to build such a test. The impossibility of building a consistent test is stated in Section \ref{Results} together with the connection with Tsirelson's equation. The proofs are relegated to Section \ref{Proofs}.

\section{Testing framework} \label{Setting}

In this section we consider an i.i.d.\ sample $(X_1,\ldots,X_n)$
of integer valued random variables with common law $\mu$. Our goal is to investigate the possibility of testing for the finiteness (or for the non finiteness) of the support of $\mu$ as $n$ goes to infinity. We refer to \cite{vva} for a detailed presentation of the notion of asymptotic statistical test.
In the following, we write $\mathbb{P}_\mu$ for the probability measure on the canonical space $\mathbb{N}^{\mathbb{N}^*}$ under which the canonical process $(X_i)_{i\ge 1}$ is a family of i.i.d.\
random variables with law $\mu$. Moreover, \noindent we denote by $\mathcal{P}_{f}$ the set of laws on $\mathbb{N}$ with finite support and by $\mathcal{P}_{\infty}$ the set of laws on $\mathbb{N}$ with infinite support. To fix idea, let us consider in this section the null hypothesis of the finiteness of the support against the alternative of the non finiteness of it (the following definitions can be adapted in an obvious way if the null and alternative hypotheses are switched).\\

\subsection{Uniform test}

In order to discriminate between the two situations, the first idea is to try to build a consistent test with asymptotic uniform level $\alpha$, with $\alpha$ given in $[0,1)$. This
means designing a rejection area $W_n$ in $\mathbb{N}^n$ such that
\begin{equation} \label{uniflevel}
\limsup_n \underset{\mu\in\mathcal{P}_f}{\text{sup}}\mathbb{P}_{\mu}\bigl[(X_1,\ldots,X_n)\in W_n\bigr] \le \alpha
\end{equation}
and for all $\mu\in\mathcal{P}_{\infty}$, 
$$\lim_n \mathbb{P}_{\mu}\bigl[(X_1,\ldots,X_n)\in W_n\bigr] =1.$$ 
This is clearly impossible. Indeed, for fixed $n$, any law $\mu_1^{\otimes n}$ with $\mu_1\in\mathcal{P}_f$ can be approached with arbitrary accuracy (in total variation norm for example) by
a law $\mu_2^{\otimes n}$ with $\mu_2\in\mathcal{P}_{\infty}$.\\

\subsection{Pointwise test}

\noindent We now consider a much weaker notion of test: we replace \eqref{uniflevel} by
$$\underset{\mu\in\mathcal{P}_{f}}{\text{sup}}\limsup_n\, \mathbb{P}_{\mu}\bigl[(X_1,\ldots,X_n)\in W_n\bigr] \le \alpha.$$
The statistical meaning of this notion of test is of course very arguable but our point of view is the following: we want to show that it is impossible to build a consistent test with given level $\al$ even in this very weak setting.\\

\begin{remark} If we do not consider all distributions with finite support but only distributions with support in
$[0,N]$ for some known $N$, a test statistics can of course be easily designed. Indeed, it is then enough to consider
the statistics $\text{max}(X_1,\ldots,X_n)\land (N+1)$: when $X_1$ follows $\mu_{1}\in\mathcal{P}_{f}$ it converges to the right endpoint of $\mu_1$ which is smaller than $N$,  when $X_1$ follows $\mu_{2}\in\mathcal{P}_{\infty}$ it converges to $N+1$.
\end{remark}

\section{Candidate test statistics for pointwise test} \label{Candidate}

We want to investigate potential test statistics for discriminating between our two situations. Since the laws in $\mathcal{P}_{f}$ and $\mathcal{P}_{\infty}$ are different because they have different support, a natural idea is to design a test statistics based on the empirical support. For example, we can split the sample of size $2n$, into two parts and compare the maximum over each subsample, that is we consider
$$S_n=\max_{1\le i\le n} X_i, \quad \tilde S_n=\max_{n+1\le i\le 2n} X_i,
\quad\text{and}\quad T_{2n}= \ind{\displaystyle \{S_n=\tilde S_n\}}.$$
If $\mu$ has finite support, it is clear that for $n$ large enough, the maxima of the two subsamples coincide and in particular
$T_{2n}$ converges in probability to $1$. 
In the infinite support case, one may expect an opposite behavior. In fact, we have the following result whose proof is given in Section \ref{Proofs}.

\begin{proposition}\label{prop1}
If $\mu$ has infinite support, then $T_{2n}$ does not converge in $\mathbb{P}_\mu$-probability to $1$.  
\end{proposition}

\noindent Therefore for any $\mu_1$ with finite support and $\mu_2$ with infinite support the statistics $T_{2n}$ shows different asymptotic behaviors under $\mu_1$ and $\mu_2$. However, this is not enough to build a consistent test since we need to investigate the asymptotic probability that $T_{2n}$ is different from $1$. In fact, the first part of Theorem \ref{prop2} in the next section implies that there exist a distribution $\mu_2$ with infinite support and a subsequence of $T_{2n}$ which converges to $1$ in probability under $\mathbb{P}_{\mu_2}$.

\section{Impossibility of testing} \label{Results}

\noindent The following theorem shows that it is not possible to build a test for the finiteness (or for the non finiteness) of the support.
Indeed, it states that if one finds a set in which any i.i.d.\ sample of size $n$ from a law with finite support $\mu_1$ is very unlikely to be, there is also one distribution $\mu_2$ with infinite support so that an  i.i.d.\ sample of it is very unlikely to be in this set, and conversely.\\

\noindent We write $\mathbb{E}_\mu$ for the expectation with respect to $\mathbb{P}_\mu$. We have the following result.

\begin{theorem}\label{prop2}
Let $\alpha\in[0,1)$, let $\phi: \mathbb{N}^*\to\mathbb{N}^*$
be an increasing map and let $A_{n}:\N^{\phi(n)}\to[0,1]$, $n\ge 1$, be a sequence of measurable functions.
\begin{enumerate}
\item
Assume that for any distribution $\mu_1$ on $\mathbb{N}$ with finite support,
$$\limsup_n\, \mathbb{E}_{\mu_1}\bigl[A_{n}(X_{1},\ldots,X_{\phi(n)}) \bigr]\leq \alpha.$$
Then, there exists a distribution 
$\mu_{2}$ on  $\mathbb{N}$ with infinite support  such that
$$\liminf_n\, \mathbb{E}_{\mu_{2}}\bigl[A_{n}(X_1,\ldots,X_{\phi(n)})\bigr]\leq \alpha.$$
\item Assume that for any distribution $\mu_2$ on $\mathbb{N}$ with infinite support,
$$\limsup_n\, \mathbb{E}_{\mu_2}\bigl[A_{n}(X_{1},\ldots,X_{\phi(n)}) \bigr]\leq \alpha.$$
Then, for all $\al'\in(\al,1)$ there exists a distribution 
$\mu_{1}$ on  $\mathbb{N}$ with finite support such that
$$\liminf_n\, \mathbb{E}_{\mu_{1}}\bigl[A_{n}(X_1,\ldots,X_{\phi(n)})\bigr]\leq \alpha'.$$
\end{enumerate}
\end{theorem}

\begin{remark}
In Theorem \ref{prop2}, we do not only consider the usual testing framework where the $A_{n}$ are indicators of a set and $\phi(n)=n$. Indeed, we allow for randomized test procedures and subsequences in the sample size. This is slightly more general and will be useful for the proof of Corollary \ref{cor1}.\\ 
\end{remark}


\noindent We now come back to Tsirelson's equation and state the result showing the impossibility of testing the hypothese that $\nu$ belongs to the union of Case 1 and Case 2 against the one that $\nu$ belongs to Case 3, and conversely.
Denote by $\mathbb{P}^*_\nu$ the law on $\T^{-\N}$ of the ``uniform solution" $(\eta^*_k)_{k\in-\N}$, given by \eqref{unif} (when $\nu_k=\nu$ for all $k$) and by $(\eta_k)_{k\in-\N}$ the canonical process on $\T^{-\N}$. We have the following corollary.

\begin{cor} \label{cor1}
Let $\alpha\in[0,1)$, let $\phi: \mathbb{N}\to\mathbb{N}$ be an increasing map
and let $B_n\subset {\T}^{\phi(n)+1}$ be a sequence of measurable sets.
\begin{enumerate}
\item
Assume that for any distribution $\nu_1$ on $\T$ belonging to Case 1 or Case 2 of Section \ref{intro},
$$\limsup_n\, \mathbb{P}^*_{\nu_1}\bigl[(\eta_0,\ldots,\eta_{-\phi(n)})\in B_{n}\bigr]\leq \alpha.$$
Then, there exists a distribution  $\nu_{2}$ belonging to Case 3 such that
$$\liminf_n\,\mathbb{P}^*_{\nu_{2}}\bigl[(\eta_0,\ldots,\eta_{-\phi(n)})\in B_{n}\bigr]\leq \alpha.$$
\item Assume that for any distribution $\nu_2$ on $\T$ belonging to Case 3 of Section \ref{intro},
$$\limsup_n\, \mathbb{P}^*_{\nu_2}\bigl[(\eta_0,\ldots,\eta_{-\phi(n)})\in B_{n}\bigr]\leq \alpha.$$
Then, for all $\al'\in(\al,1)$ there exists a distribution $\nu_1$ belonging to Case 1 or Case 2 of Section \ref{intro} such that
$$\liminf_n\,\mathbb{P}^*_{\nu_{1}}\bigl[(\eta_0,\ldots,\eta_{-\phi(n)})\in B_{n}\bigr]\leq \alpha'.$$
\end{enumerate}
\end{cor}

\section{Proofs}\label{Proofs}

\subsection{Proof of Proposition \ref{prop1}}

Let $\mu$ be a probability measure on $\N$ and assume that $\mathbb{P}_\mu[T_{2n}=1]\to 1$. Since $S_n$ and $\tilde S_n$ are independent with the same law,
it implies that there exists a sequence of integers $k_n$ such that
$\mathbb{P}_\mu[S_n=k_n]\to 1$.\\
\noindent Let us show that $k_n$ is a bounded sequence, which will imply that the support of $\mu$ is finite. If $k_n$ is not bounded, then there exist a subsequence $n_\ell$ such that
$k_{1+n_\ell}> k_{n_\ell}$ for all $\ell$ and $k_{n_\ell}\to\infty$. Therefore $$\mathbb{P}_\mu[S_{n_\ell}=k_{n_\ell}] \to 1$$
and
$$ \mathbb{P}_\mu[S_{1+n_\ell}=k_{n_\ell}] \to 0$$
because $P\bigpar{S_{1+n_\ell}=k_{1+n_\ell}} \to 1$ and $k_{1+n_\ell}\not=k_{n_\ell}$.
Moreover, we have that
$$\mathbb{P}_\mu[S_{1+n_\ell}=k_{n_\ell}]-\mathbb{P}_\mu[S_{n_\ell}=k_{n_\ell}]$$ is equal to 
$$\mathbb{P}_\mu[S_{n_\ell}=k_{n_\ell}] \mu([0,k_{n_\ell}]) + \mathbb{P}_\mu[S_{n_\ell}<k_{n_\ell}] \mu(k_{n_\ell})
- \mathbb{P}_\mu[S_{n_\ell}=k_{n_\ell}].$$
This absolute value of this last quantity is smaller than $\mu([k_{n_\ell},+\infty))$ which goes to zero as $n$ goes to infinity. This
shows the contradiction.

\subsection{Proof of Theorem \ref{prop2}}

\subsubsection{Proof of Part 1 in Theorem \ref{prop2}}

We recursively define a sequence of distributions $\mu_{2}^n$, $n\ge 1$, and a sequence of integers $\psi(n)$, $n\ge 0$.\\
\noindent $-$ At rank $n=1$, we consider $\mu_{2}^1$ the Dirac measure at point $0$, $\psi(0)=1$ and $\psi(1)=1$.\\
\noindent $-$ At rank $n>1$, we define $\mu_{2}^n$ and $\psi(n)$. The law $\mu_{2}^n$ is the distribution with discrete support $\{0,\ldots,n-1\}$ defined by $$\mu^n_2(k)=\frac{c_n}{(\phi\circ\psi(k))^2},\quad k\in\{0,\ldots,n-1\},$$
with $c_n$ such that $$\sum_{k=0}^{n-1}\frac{c_n}{(\phi\circ\psi(k))^2}=1.$$ The integer $\psi(n)$ is taken so that $$\psi(n)>\text{max}(\psi(n-1),n^2)$$ and for all $m\geq\psi(n)$,
$$\mathbb{E}_{\mu^n_2}\bigl[A_m(X_1,\ldots,X_{\phi(m)})\bigr]\leq \alpha+\frac{1}{n}.$$
Note that finding such a $\psi(n)$ is always possible thanks to the assumption on the sequence $(A_m)_m$. In particular, the sequence $\psi(n)$ is increasing and satisfies
$$\sum_{k=0}^{+\infty}\frac{1}{(\phi\circ\psi(k))^2}<\infty.$$

\noindent Now define the distribution $\mu_{2}$ on $\mathbb{N}$ with infinite support by
$$\mu_2(k)=\frac{c}{(\phi\circ\psi(k))^2},~k\in\mathbb{N},$$
with $c$ such that $$\sum_{k=0}^{+\infty}\frac{c}{(\phi\circ\psi(k))^2}=1.$$
We have that
$$\mathbb{E}_{\mu_2}\bigl[A_{\psi(n)}(X_1,\ldots,X_{\phi\circ\psi(n)})\bigr]$$
is smaller than
\begin{align*}
&\mathbb{E}_{\mu_2}\bigl[A_{\psi(n)}(X_1,\ldots,X_{\phi\circ\psi(n)}) \mid \bigcap_{i=1}^{\phi\circ\psi(n)}\set{X_i\le n-1}\bigr] +
\mathbb{P}_{\mu_2}\bigl[ \bigcup_{i=1}^{\phi\circ\psi(n)}\set{X_i\ge n} \bigr]\\
&=\mathbb{E}_{\mu^n_2}\bigl[A_{\psi(n)}(X_1,\ldots,X_{\phi\circ\psi(n)})\bigr] +
\mathbb{P}_{\mu_2}\Bigl[ \bigcup_{i=1}^{\phi\circ\psi(n)}\set{X_i\ge n} \Bigr]\\
& \le \alpha+\frac{1}{n}+\phi\circ\psi(n)\sum_{k=n}^{+\infty}\frac{1}{(\phi\circ\psi(k))^2}.
\end{align*}
Using the fact that $\psi(n)$ is increasing and the inequality $\psi(n)>n^2$, we obtain
$$\alpha+\frac{1}{n}+\phi\circ\psi(n)\sum_{k=n}^{+\infty}\frac{1}{(\phi\circ\psi(k))^2}\leq\alpha+\frac{1}{n}+\sum_{k=n}^{+\infty}\frac{1}{\phi\circ\psi(k)}\leq \alpha+\frac{1}{n}+\sum_{k=n}^{+\infty}\frac{1}{k^2}.$$
 This quantity goes to $\alpha$ as $n$ goes to infinity, which gives the result.
 
\subsubsection{Proof of Part 2 in Theorem \ref{prop2}}
Let $\al'\in(\al,1)$.
Assume that for any distribution 
$\mu_{1}$ on  $\mathbb{N}$ with finite support
$$\liminf_n\, \mathbb{E}_{\mu_{1}}\bigl[A_{n}(X_1,\ldots,X_{\phi(n)})\bigr]> \alpha'.$$
This last inequality is equivalent to
$$\limsup_n\, \mathbb{E}_{\mu_{1}}\bigl[1-A_{n}(X_1,\ldots,X_{\phi(n)})\bigr]< 1-\alpha'.$$
According to point 1.\ of Theorem \ref{prop2}, there exists a distribution $\mu_2$ on $\N$ with infinite
support such that
$$\liminf_n\, \mathbb{E}_{\mu_{2}}\bigl[1-A_{n}(X_1,\ldots,X_{\phi(n)})\bigr]\le 1-\alpha',$$
that is
$$\limsup_n\, \mathbb{E}_{\mu_{2}}\bigl[A_{n}(X_1,\ldots,X_{\phi(n)})\bigr]\ge \alpha'$$
wich gives the contradiction since $\al'>\al$.

\subsection{Proof of Corollary \ref{cor1}}
Let $f$ be an injection from $\N$ into $\T\cap \Q/\Z$. If $\mu$ is a probability measure on $\N$, denote
by $f(\mu)=\mu\circ f^{-1}$ the image of $\mu$ by $f$.
On the one hand, one has:\\
- if $\mu$ has finite support, $f(\mu)$ belongs to Case 1 or Case 2.\\
- if $\mu$ has infinite support, $f(\mu)$ belongs to Case 3 since $f(\mu)$ has infinite support.\\
On the other hand, because of the definition of $\mathbb{P}^*_\nu$, one has
$$\mathbb{P}^*_{f(\mu)}\bigl[(\eta_0,\ldots,\eta_{-\phi(n)})\in B_{n}\bigr]
=\mathbb{E}_{\mu}\bigl[A_n(X_1,\dots,X_{\phi(n)})\bigr]$$
where $A_n(X_1,\dots,X_{\phi(n)})$ is equal to
$$\int_0^1 \ind{B_n}\bigl(u,u-f(X_1), u-f(X_1)-f(X_2),\dots, u-f(X_1)-f(X_2)-\dots
-f(X_{\phi(n)})\bigr)\,du.
$$
Using these two facts, together with Theorem \ref{prop2}, the proof of Corollary \ref{cor1} is easily completed.
\section*{Acknowledgments}
We are very grateful to Kouji Yano and Marc Yor for introducing us to the statistical question linked to Tsirelson's equation.

\end{document}